\documentclass[a4paper,12pt]{article}
\usepackage[T2A]{fontenc}
\usepackage{amssymb,amsmath,amsfonts,amsthm}
\usepackage[utf8]{inputenc} 
\usepackage[english]{babel}
\usepackage[pdftex,unicode]{hyperref}
\usepackage{indentfirst} 
  
\newtheorem{rem}{Remark}
\newtheorem{theorem}{Theorem}
\newtheorem{corollary}{Corollary}

\usepackage{longtable}
\usepackage{multicol}
\author{Kostousov, K.V.}
\title{Symmetrical 2-extensions of the 3-dimensional grid. I}

\date{\vspace{-5ex}}

\begin{document}

\newcommand{\Addresses}{{
		\bigskip
		\footnotesize
		
		Kirill Viktorovich Kostousov, \textsc{Krasovskii Institute of Mathematics and Mechanics, S.Kovalevskaya Str, 16, 620108, Yekaterinburg, Russia;}\par\nopagebreak
		\textit{E-mail address: } \texttt{kkostousov@gmail.com}
}}

\maketitle	

\begin{abstract}
For a positive integer $d$, a connected graph $\Gamma$ is a symmetrical
2-extension of the $d$-dimensional grid
$\Lambda^d$ if there exists a vertex-tran\-sitive group $G$ of 
automorphisms
of $\Gamma$ and its
imprimitivity system $\sigma$ with blocks of order 2
such that there exists an isomorphism $\varphi$ of the
quotient graph $\Gamma/\sigma$ onto $\Lambda^d$.
The tuple $(\Gamma, G, \sigma, \varphi)$ with specified components is 
called a realization
of the symmetrical 2-extension $\Gamma$ of the grid $\Lambda^{d}$.
Two realizations $(\Gamma_1, G_1,$ $\sigma_1, \varphi_1)$ and 
$(\Gamma_2, G_2, \sigma_2, \varphi_2)$
are called equivalent
if there exists an isomorphism of the graph $\Gamma_1$ onto $\Gamma_2$ 
which maps $\sigma_1$ onto $\sigma_2$.
In \cite{IV3}, V. Trofimov proved that, up to equivalence, there are 
only finitely many realizations of
symmetrical $2$-extensions of $\Lambda^{d}$
for each positive integer $d$. In \cite{dim2partI} and 
\cite{dim2partII}, E. Konovalchik and K. Kostousov
found all, up to equivalence,
realizations of symmetrical 2-extensions of the grid $\Lambda^2$.
In this work we found all, up to equivalence, 
realizations $(\Gamma, G, \sigma, \varphi)$ of symmetrical 2-extensions
of the grid $\Lambda^3$ for which only the trivial automorphism of 
$\Gamma$
preserves all blocks of $\sigma$ (we prove that there are 5573 such 
realizations,
and that among corresponding graphs $\Gamma$ there are 5350 pairwise 
non-isomorphic).
Our work was performed using «Uran» supercomputer of IMM UB RAS.


\end{abstract}

\section{Introduction}
Recall that for a positive integer $d$ a $d$-dimensional grid $\Lambda^{d}$ is a graph,
which vertexes are iteger tuples $(a_{1},\ldots a_{d})$ and two vertexes
$(a_{1}',\ldots,a_{d}')$ and  $(a_{1}'',\ldots,a_{d}'')$ are adjacent if and only if 
$|a_{1}'-a_{1}''|+\ldots+|a_{d}'-a_{d}''|=1$.
For a finite graph $\Delta$ define a connected graph $\Gamma$ to be {\it 
a symmetrical extension of $\Lambda^{d}$ by $\Delta$} if there exists a vertex-transitive
group $G$ of automorphisms of $\Gamma$ and an imprimitivity system $\sigma$ of $G$ on $V(\Gamma)$
such that subgraphs of $\Gamma$, generated by blocks of $\sigma$ are isomorphic to $\Delta$
and there exists an isomorphism of $\Gamma/\sigma$ (i.e. of factor-graph of $\Gamma$ by partition $\sigma$ of its vertex set)
onto $\Lambda^{d}$. A tuple $(\Gamma, G, \sigma, \varphi)$ with specified components if called a {\it
realizaion of symmetrical extension $\Gamma$ of the grid $\Lambda^{d}$ by the graph $\Delta$}.
For a positive integer $q$ a graph $\Gamma$ is called {\it a symmetrical $q$-extension of the grid
 $\Lambda^{d}$}, if $\Gamma$ is a symmetrical extension of the grid $\Lambda^d$  by some graph $\Delta$, such
that $|V(\Delta)|=q$. In this situation the tuple $(\Gamma, G, \sigma, \varphi)$ with specified components
is called {\it a realization of the symmetrical $q$-extension} $\Gamma$ of the grid $\Lambda^{d}$,
and $\Gamma$ we call a graph of this realization.
Along with purely mathematical interest,
symmetrical $q$-extensions of the grid $\Lambda^{d}$ for small $d\geq 1$ and $q>1$ are iteresting
for molecular crystallography and some physical theories (see \cite{IV2}).
For the crystallography symmetrical $2$-extensions of grids $\Lambda^d$
are of the most interest. They naturally arise when considering ``molecular'' crystals,
whose ``molecules'' consist of two ``atoms'' or, more generally, have a accentuated axis.

It is natural to consider realizations of symmetric $q$-extensions
of the grid $\Lambda^{d}$ up to equivalence defined as follows (see \cite {IV2}).
We call two such realizations $R_1=(\Gamma_1, G_1,$ $\sigma_1, \varphi_1)$ and $R_2=(\Gamma_2, G_2, \sigma_2, \varphi_2)$ 
{\it equivalent}, and we will write $R_1 \sim R_2$ if there exists an isomorphism of the graph 
$\Gamma_1$ to the graph $\Gamma_2$ which maps $\sigma_1$ onto $\sigma_2$.
The realizaion $(\Gamma, G, \ sigma, \ varphi)$ of the symmetrical $q$-extension of the grid $\Lambda^{d}$
we call {\it maximal} if $G=\mathrm{Aut}_\sigma(\Gamma)$ is the group of all automorphisms of the graph 
$\Gamma$ which preserve the partition $\sigma$.
It is clear that each realization of the symmetrical $q$-extension of the grid $\Lambda^{d}$ 
has an equivalent maximal realization (unique up to equivalence).
V.I. Trofimov proved that for an arbitrary positive integer $d$
up to equivalence, there is only a finite number of realizations of symmetrical 2-extensions
of $d$-dimensional grid (see \cite[Theorem 2]{IV3}).
An algorithm for constructing of these extensions is also proposed in \cite{IV3}.

Using this algorithm, in \cite{dim2partI} and \cite{dim2partII} was found all, up to equivalence,
realizations of symmetrical 2-extensions of the grid $\Lambda^2 $ (162 realizations).
Among the graphs of these realizations, there are exactly 152 pairwise nonisomorphic graphs.

For an arbitrary realization $(\Gamma,  G, \sigma, \varphi)$ of the symmetrical 2-extension of the grid $\Lambda^{d}$
and an arbitrary pair of adjacent vertices $B_1, B_2$ of the graph $\Gamma/\sigma$
the set of edges of the graph $\Gamma$, one end of which lies in $B_1$, and the other in $B_2$,
we will call a {\it connection}.
The following types of connections are possible:
{\it type $1$} --- a single edge;
{\it type $2_{||}$} --- two non-adjacent edges;
{\it type $2_V$} --- two adjacent edges;
{\it type $3$} --- three edges;
{\it type $4$} --- full connection (4 edges).
A realization, which necessarily have connections of types not equal to  $2_{||}$ and $4$, we will call by {\it realization of class} I.
A realization, which have connections only of types $2_{||}$ and $4$ (maybe only of one type), we will call by {\it realizations of class} II.
By Proposition 4 from \cite{IV3} realizations of class I,
are exactly the realizations of symmetrical 2-extensions of the grid $\Lambda^d $ such that
only a trivial automorphism of their graph fixes all blocks (as whole).

All 162 realizations of symmetrical 2-extensions of the grid $\Lambda^2$
are distributed in classes I and II as follows: 87 realizations of class I (see \cite {dim2partI})
and 75 realizations of class II (see \cite {dim2partII}).
This paper is devoted to the description of all, up to equivalence,
realizations of symmetrical 2-extensions of the grid $\Lambda^3$ of class I.

A realization of symmetrical extension of the grid $\Lambda^{d}$ by the graph $K_2$ (full graph on two vertices)
we will call the {\ it saturated} realization of the symmetrical 2-extension of the grid $\Lambda^{d}$.
Accordingly, a realization of symmetrical extension of the grid $\Lambda^{d}$ by the graph, complemented to $K_2$
we will call {\ it non-saturated} realization of the symmetrical 2-extension of the grid $\Lambda^{d}$.

In this paper, we have shown that, up to equivalence,
there are 5573 realizations of symmetrical 2-extensions of the grid $\Lambda^{3}$ of
class I, among which 2872 are saturated and 2701 are non-saturated (see Theorem 1 and Corollary 1).
Among the graphs of saturated realizations of symmetrical 2-extensions of the grid $\Lambda^{3}$ of class I 
there are exactly 2792 pairwise nonisomorphic;
among the graphs of non-saturated realizations of class I there are
exactly 2594 pairwise nonisomorphic;
and among all graphs of realizations of class I there are
5350 pairwise nonisomorphic (see Corollary 2).

In Sec. \ref{main_res} we give the desciption of all up to equivalence realizations of symmetrical 2-extensions
of the grid $\Lambda^{3}$ of class I (Theorem 1 and Corollary 1). It is obtained in this article
 using the approach from \cite{IV3} implemented in GAP\cite{gap} (Algorithms 1 and 2 from \cite {dim2partI}). 
Sec. \ref{auxiliary} contains preliminary results.

Our work was performed using «Uran» supercomputer of IMM UB RAS.

\section{Preliminaries}\label{auxiliary}

Each vertex-transitive group of automorphisms of
$\Lambda^{3}$ is generated by the stabilizer in this 
group of the vertex $(0,0,0)$ and six elements of this group,
translating the vertex $(0,0,0) $ to the vertices adjacent to it.
Based on this, using GAP \cite{gap}, we have listed
all conjugate classes of vertex-transitive subgrous of the group $\mathrm{Aut}(\Lambda^{3})$.
It turned out that there are 786 such classes 
whose system of representatives we denote by $\textbf{H}=\big\{H_{1},\ldots, H_{786}\big\}$.
The groups $H_{1},\ldots, H_{786}$ are given in the Table 2 below by their generating systems.
The following notation is used for certain automorphisms of the grid $\Lambda^{3}$:
\begin{center}
\begin{tabular}{cc}
$r_x: (x,y,z)\mapsto (x,-z,y),$ &
$r_y: (x,y,z)\mapsto (z,y,-x),$ \\
$r_z: (x,y,z)\mapsto (y,-x,z),$ &
$m_x: (x,y,z)\mapsto (-x,y,z),$ \\
$m_y: (x,y,z)\mapsto (x,-y,z),$ &
$m_z: (x,y,z)\mapsto (x,y,-z),$ \\
$i: (x,y,z)\mapsto (-x,-y,-z),$ & 
$t_x: (x,y,z)\mapsto (x+1,y,z),$ \\
$t_y: (x,y,z)\mapsto (x,y+1,z),$ &
$t_z: (x,y,z)\mapsto (x,y,z+1),$
\end{tabular}
\end{center}
где $x,y,z\in\mathbb Z$.
\smallskip

\begin{rem} 
In the natural embedding of the grid $\Lambda^{3}$ in the Euclidean space
$\mathbb R^3$ each automorphism
$g\in \mathrm{Aut}(\Lambda^{3})$ is induced by the only isometry $\tilde{g} $ of this space.
The isometries that induce the above automorphisms of $\Lambda^{3}$ have the following geometric meaning:
$\tilde{r_x}, \tilde{r_y}, \tilde{r_z}$ are rotations by the angle $\frac{\pi}{2}$ around coordinate axes $x$, $y$ and $z$ correspondingly,
$\tilde{m_x}, \tilde{m_y}, \tilde{m_z}$ are reflections relative coordinate planes,
$\tilde i$ --- central symmetry abount the origin,
$\tilde{t_x}, \tilde{t_y}, \tilde{t_z}$ are translations by 1 along axes $x$, $y$ and $z$ correspondingly.\smallskip
\end{rem}

Using GAP, we constructed and compared 786 stabilizers \linebreak $\{H_{(0,0,0)} : H\in \textbf{H}\}$ for conjugation in $\mathrm{Aut}(\Lambda^{3})_{(0,0,0)}$. 
It turned out that up to conjugation there are only 33 such stabilizers.
We give them in Table 1 by their generators (column 3).  
For each of the 33 groups, the abstract group structure is given (column 2).

\hspace*{97mm} \mbox{T\ a\ b\ l\ e\ \  1}\\
\centerline{\small\bf
  Stabilizers of vertex $(0,0,0)$ in vertex-transitive subgroups of $\mathrm{Aut}(\Lambda^3)$}
  \centerline{\small\bf up to conjugation in $\mathrm{Aut}(\Lambda^3)_{(0,0,0)}$}
\begin{longtable}{llll}
\textnumero & Group structure & Generators  \\
\hline\endhead
1& $1$ & 1\\
2& $C_2$ & $\langle i \rangle$\\
3& $C_2$ & $\langle m_z \rangle$\\
4& $C_2$ & $\langle r_z^2 r_x \rangle$\\
5& $C_2$ & $\langle r_z^2 \rangle$\\
6& $C_2$ & $\langle m_z r_x \rangle$\\
7& $C_3$ & $\langle r_y^{-1} r_z^{-1} \rangle$\\
8& $C_2\times C_2$ & $\langle r_y^2, r_z^2 \rangle$\\
9& $C_2\times C_2$ & $\langle i, r_z^2 \rangle$\\
10& $C_2\times C_2$ & $\langle m_z r_x^{-1}, r_x^2 \rangle$\\
11& $C_2\times C_2$ & $\langle m_x, r_y^2 r_x \rangle$\\
12& $C_2\times C_2$ & $\langle i, r_z^2 r_x \rangle$\\
13& $C_2\times C_2$ & $\langle r_y^2 r_x, r_z^2 r_x \rangle$\\
14& $C_2\times C_2$ & $\langle m_x, r_z^2 \rangle$\\
15& $C_4$ & $\langle r_x^{-1}, r_x^2 \rangle$\\
16& $C_4$ & $\langle m_x r_x^{-1} \rangle$\\
17& $C_6$ & $\langle i, m_x r_y^{-1} r_x^{-1} \rangle$\\
18& $S_3$ & $\langle m_x r_y, r_y^{-1} r_z^{-1} \rangle$\\
19& $S_3$ & $\langle r_y^2 r_z, r_z^2 r_x \rangle$\\
20& $C_2\times C_2\times C_2$ & $\langle i, r_y^2 r_x, r_z^2 r_x \rangle$\\
21& $C_2\times C_2\times C_2$ & $\langle i, r_y^2, r_z^2 \rangle$\\
22& $C_4\times C_2$ & $\langle i, m_x r_x^{-1} \rangle$\\
23& $D_8$ & $\langle m_x r_x^{-1}, r_z^2 r_x \rangle$\\
24& $D_8$ & $\langle r_z^2, r_z^2 r_x \rangle$\\
25& $D_8$ & $\langle m_y, m_z r_x^{-1}, r_x^2 \rangle$\\
26& $D_8$ & $\langle m_x r_x^{-1}, r_z^2 \rangle$\\
27& $A_4$ & $\langle r_y r_x^{-1}, r_y^2, r_z^2 \rangle$\\
28& $D_{12}$ & $\langle i, r_y^2 r_z, r_z^2 r_x \rangle$\\
29& $C_2\times D_8$ & $\langle i, r_z^2, r_z^2 r_x \rangle$\\
30& $C_2\times A_4$ & $\langle i, m_x r_y^{-1} r_x^{-1}, r_y^2, r_z^2 \rangle$\\
31& $S_4$ & $\langle m_x r_x^{-1}, m_x r_z, r_z^2 \rangle$\\
32& $S_4$ & $\langle r_y^2 r_z, r_z^2, r_z^2 r_x \rangle$\\
33& $C_2\times S_4$ & $\langle i, r_y^2 r_z, r_z^2, r_z^2 r_x \rangle$\\
\end{longtable} 

Each group $H\in\textbf{H}$ is identified with some space group,
and, therefore, has a point group $\mathrm{P}(H)$ and a translation basis (see \cite{cryst}).
Using GAP, we verified that the set of point groups $\{\mathrm{P}(H) : H\in \textbf{H}\}$
up to conjugation in $\mathrm{P}(\mathrm{Aut}(\Lambda^{3}))$ is equal to the set of 33 stabilizers, 
given in Table 1.
In column 1 of Table 2 below, we give the set of groups $\textbf{H}$ defined by their generators.
In column 2 for each group $H\in \textbf{H}$ we give the \textnumero \, of group from Table 1 conjugate to
the stabilizer $H_{(0,0,0)}$ in $\mathrm{Aut}(\Lambda^{3})_{(0,0,0)}$.
In column 3 for each group $H\in \textbf{H}$ we give the \textnumero \, of group from Table 1, conjugate to the point group 
$\mathrm{P}(H)$ in $\mathrm{P}(\mathrm{Aut}(\Lambda^{3}))$.
In column 4 for each group $H\in\textbf{H}$ we give its translation basis.
The groups in Table 2 are sorted lexicographically first by \textnumero \, in column 2 and then by \textnumero \, in column 3.

\hspace*{97mm} \mbox{T\ a\ b\ l\ e\ \  2}\\
\centerline{\small\bf
 Conjugate classes of vertex-transitive}
  \centerline{\small\bf subgroups of the group $\mathrm{Aut}(\Lambda^3)$ representatives}
{\tiny
\begin{longtable}{llll}
$H$ & $H_{0,0,0}$ & $P(H)$ & Translation basis of $H$\\
\hline\endhead
$H_{1}=\langle t_x, t_y^{-1}, t_z^{-1} \rangle$ & 1 & 1 & [ [ 1, 0, 0 ], [ 0, 1, 0 ], [ 0, 0, 1 ] ]\\
$H_{2}=\langle i t_z, i t_z^{-1}, t_x, t_y^{-1} \rangle$ & 1 & 2 & [ [ 1, 0, 0 ], [ 0, 1, 0 ], [ 0, 0, 2 ] ]\\
$H_{3}=\langle i t_y, i t_y^{-1}, i t_z, i t_z^{-1}, t_x \rangle$ & 1 & 2 & [ [ 1, 0, 0 ], [ 0, 1, 1 ], [ 0, 0, 2 ] 
 ]\\
$H_{4}=\langle i t_x, i t_x^{-1}, i t_y, i t_y^{-1}, i t_z, i t_z^{-1} \rangle$ & 1 & 2 & 
[ [ 1, 0, 1 ], [ 0, 1, 1 ], [ 0, 0, 2 ] ]\\
$H_{5}=\langle m_x t_y^{-1}, t_x, t_z^{-1} \rangle$ & 1 & 3 & [ [ 1, 0, 0 ], [ 0, 2, 0 ], [ 0, 0, 1 ] ]\\
$H_{6}=\langle m_y t_y, m_y t_y^{-1}, t_x, t_z^{-1} \rangle$ & 1 & 3 & [ [ 1, 0, 0 ], [ 0, 2, 0 ], [ 0, 0, 1 ] ]\\
$H_{7}=\langle m_y t_y, m_y t_y^{-1}, m_y t_z^{-1}, t_x \rangle$ & 1 & 3 & [ [ 1, 0, 0 ], [ 0, 1, 1 ], [ 0, 0, 2 ] ]\\
$H_{8}=\langle m_y t_x, m_y t_z^{-1}, t_y^{-1} \rangle$ & 1 & 3 & [ [ 1, 0, 1 ], [ 0, 1, 0 ], [ 0, 0, 2 ] ]\\
\end{longtable} }

\section{Main result}\label{main_res}

We have done a computer implementation of the
proposed in \cite{IV3} approach, which can be called a coordinatization
of symmetrical extensions of graphs.
According to it, the realization of symmetrical 2-extension of the grid
$\Lambda^{3}$ of class I can be defined by a triple $H, L, X$, where $H$ is a vertex-transitive
subgroup of $\mathrm {Aut}(\Lambda^{3})$, $L$ is subgroup of index $2$ of the stabilizer
of the vertex $(0,0,0)$ in $H$, and $X$ is some subset of elements of $H$ mapping
the vertex $(0,0,0)$ of $\Lambda^{3}$ to some its adjacent vertexes (for details, see \cite{dim2partI}).

Algorithm 1 (Generating of all saturated realizations of symmetrical 2-extensions of $\Lambda^{2}$)
and Algorithm 2 (Check two realizations for equivalence) from \cite{dim2partI},
we adapted (without essential changes) and applied to $\Lambda^{3}$.
A list of saturated realizations generated by Algorithm 1 and thinned out by Algorithm 2,
contains 2872 realizations (given in Table 4 below). 
When thinning in each class of equivalent realizations the
realization with maximal by inclusion group $H_i$ was selected.
Due to this, the realizations in the resulting list are maximal.

We split the set of all realizations of the symmetrical 2-extensions of $\Lambda^{3}$
of the class I into subclasses, defined by the types of connections in the neighbourhood of vertex.
In the first column of Table 3 below, we give all occurring combinations of connection types in the neighbourhood of vertex (59 combinations).
Combinations are of the form $x_1 x_2\_ y_1 y_2\_ z_1 z_2$, where 
$x_1$ is the type of the first connection by the first direction
(the grid $\Lambda^{3}$, and therefore a 2-extension, has three directions along coordinate axes),
$x_2$ is the type of the second connection by the first direction,
$y_1$ is the type of the first connection by the second direction,
$y_2$ is the type of the second connection by the second direction,
$z_1$ is the type of the first connection by the third direction,
$z_2$ is the type of the second connection by the third direction.
Here, for each extension, the numbering of directions and connections within the direction is performed so that
the combination turned out to be lexicographically minimal.
At the pictures of combinations in the first column of Table 3
the first direction is shown horizontally
(first connection to the left, second to the right),
the second direction is shown vertically
(bottom connection is first, top connection is second),
the third direction is shown diagonally
(Bottom-left connection is first, top-right connection is second).
For each combination of connection types, the remaining columns in Table 3 contain
pictures of all found corresponding extensions of vertex neighborhood up to equivalence.
At these pictures, the edges in blocks are not shown because
we use these types of vertex neighborhood both for saturated and non-saturated realizations.

\hspace*{97mm} \mbox{T\ a\ b\ l\ e\ \  3}\\
\centerline{\bf Vertex neighbourhood extensions for}\\
 \centerline{\small\bf 2-extensions of $\Lambda^{3}$ of class I}
\begin{longtable}{llll}
connection types & \multicolumn{3}{c}{vertex neighbourhood extensions}  \\
\hline\endhead
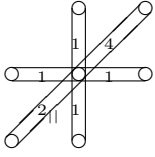
\begin{tabular}{c} \input{pictures/det1.tex} \\$11\_11\_11$\end{tabular}& \begin{tabular}{c} \input{pictures/ball1_1.tex} \\1A\end{tabular}& \begin{tabular}{c} \input{pictures/ball1_2.tex} \\1B\end{tabular}\\
\begin{tabular}{c} \input{pictures/det2.tex} \\$11\_11\_2_{||}2_{||}$\end{tabular}& \begin{tabular}{c} \input{pictures/ball2_1.tex} \\2A\end{tabular}& \begin{tabular}{c} \input{pictures/ball2_2.tex} \\2B\end{tabular}\\
\begin{tabular}{c} \input{pictures/det3.tex} \\$11\_11\_2_{||}4$\end{tabular}& \begin{tabular}{c} \input{pictures/ball3_1.tex} \\3A\end{tabular}& \begin{tabular}{c} \input{pictures/ball3_2.tex} \\3B\end{tabular}\\
\begin{tabular}{c} \input{pictures/det4.tex} \\$11\_11\_33$\end{tabular}& \begin{tabular}{c} \input{pictures/ball4_1.tex} \\4A\end{tabular}& \begin{tabular}{c} \input{pictures/ball4_2.tex} \\4B\end{tabular}\\
\begin{tabular}{c} \input{pictures/det5.tex} \\$11\_11\_44$\end{tabular}& \begin{tabular}{c} \input{pictures/ball5_1.tex} \\5A\end{tabular}& \begin{tabular}{c} \input{pictures/ball5_2.tex} \\5B\end{tabular}\\
\begin{tabular}{c} \input{pictures/det6.tex} \\$11\_12_{||}\_12_{||}$\end{tabular}& \begin{tabular}{c} \input{pictures/ball6_1.tex} \\6\end{tabular}\\
\begin{tabular}{c} \input{pictures/det7.tex} \\$11\_13\_13$\end{tabular}& \begin{tabular}{c} \input{pictures/ball7_1.tex} \\7A\end{tabular}& \begin{tabular}{c} \input{pictures/ball7_2.tex} \\7B\end{tabular}\\
\begin{tabular}{c} \input{pictures/det8.tex} \\$11\_14\_14$\end{tabular}& \begin{tabular}{c} \input{pictures/ball8_1.tex} \\8\end{tabular}\\
\begin{tabular}{c} \input{pictures/det9.tex} \\$11\_2_{||}2_{||}\_2_{||}2_{||}$\end{tabular}& \begin{tabular}{c} \input{pictures/ball9_1.tex} \\9\end{tabular}\\
\begin{tabular}{c} \input{pictures/det10.tex} \\$11\_2_{||}2_{||}\_2_{||}4$\end{tabular}& \begin{tabular}{c} \input{pictures/ball10_1.tex} \\10\end{tabular}\\
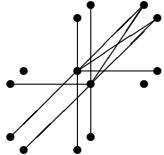
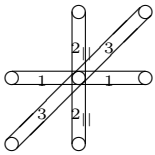
\begin{tabular}{c} \input{pictures/det11.tex} \\$11\_2_{||}2_{||}\_33$\end{tabular}& \begin{tabular}{c} \input{pictures/ball11_1.tex} \\11\end{tabular}\\
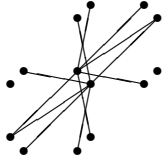
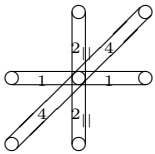
\begin{tabular}{c} \input{pictures/det12.tex} \\$11\_2_{||}2_{||}\_44$\end{tabular}& \begin{tabular}{c} \input{pictures/ball12_1.tex} \\12\end{tabular}\\
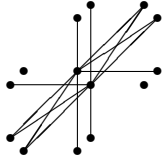
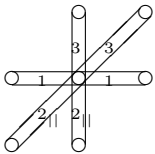
\begin{tabular}{c} \input{pictures/det13.tex} \\$11\_2_{||}3\_2_{||}3$\end{tabular}& \begin{tabular}{c} \input{pictures/ball13_1.tex} \\13\end{tabular}\\
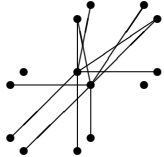
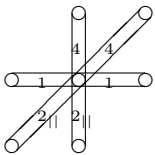
\begin{tabular}{c} \input{pictures/det14.tex} \\$11\_2_{||}4\_2_{||}4$\end{tabular}& \begin{tabular}{c} \input{pictures/ball14_1.tex} \\14\end{tabular}\\
\begin{tabular}{c} \input{pictures/det15.tex} \\$11\_2_{||}4\_33$\end{tabular}& \begin{tabular}{c} \input{pictures/ball15_1.tex} \\15\end{tabular}\\
\begin{tabular}{c} \input{pictures/det16.tex} \\$11\_2_{||}4\_44$\end{tabular}& \begin{tabular}{c} \input{pictures/ball16_1.tex} \\16\end{tabular}\\
\begin{tabular}{c} \input{pictures/det17.tex} \\$11\_2_V2_V\_2_V2_V$\end{tabular}& \begin{tabular}{c} \input{pictures/ball17_1.tex} \\17A\end{tabular}& \begin{tabular}{c} \input{pictures/ball17_2.tex} \\17B\end{tabular}\\
\begin{tabular}{c} \input{pictures/det18.tex} \\$11\_33\_33$\end{tabular}& \begin{tabular}{c} \input{pictures/ball18_1.tex} \\18A\end{tabular}& \begin{tabular}{c} \input{pictures/ball18_2.tex} \\18B\end{tabular}\\
\begin{tabular}{c} \input{pictures/det19.tex} \\$11\_33\_44$\end{tabular}& \begin{tabular}{c} \input{pictures/ball19_1.tex} \\19\end{tabular}\\
\begin{tabular}{c} \input{pictures/det20.tex} \\$11\_34\_34$\end{tabular}& \begin{tabular}{c} \input{pictures/ball20_1.tex} \\20\end{tabular}\\
\begin{tabular}{c} \input{pictures/det21.tex} \\$11\_44\_44$\end{tabular}& \begin{tabular}{c} \input{pictures/ball21_1.tex} \\21\end{tabular}\\
\begin{tabular}{c} \input{pictures/det22.tex} \\$12_{||}\_12_{||}\_2_{||}2_{||}$\end{tabular}& \begin{tabular}{c} \input{pictures/ball22_1.tex} \\22\end{tabular}\\
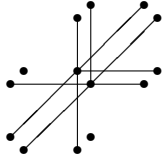
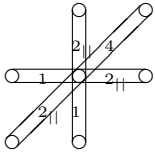
\begin{tabular}{c} \input{pictures/det23.tex} \\$12_{||}\_12_{||}\_2_{||}4$\end{tabular}& \begin{tabular}{c} \input{pictures/ball23_1.tex} \\23\end{tabular}\\
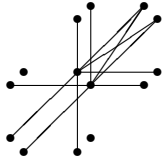
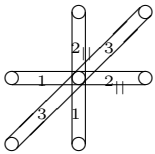
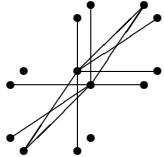
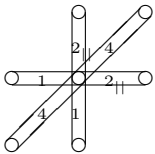
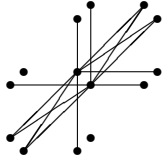
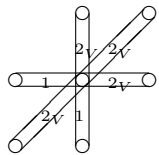
\begin{tabular}{c} \input{pictures/det24.tex} \\$12_{||}\_12_{||}\_33$\end{tabular}& \begin{tabular}{c} \input{pictures/ball24_1.tex} \\24\end{tabular}\\
\begin{tabular}{c} \input{pictures/det25.tex} \\$12_{||}\_12_{||}\_44$\end{tabular}& \begin{tabular}{c} \input{pictures/ball25_1.tex} \\25\end{tabular}\\
\begin{tabular}{c} \input{pictures/det26.tex} \\$12_V\_12_V\_2_V2_V$\end{tabular}& \begin{tabular}{c} \input{pictures/ball26_1.tex} \\26A\end{tabular}& \begin{tabular}{c} \input{pictures/ball26_2.tex} \\26B\end{tabular}& \begin{tabular}{c} \input{pictures/ball26_3.tex} \\26C\end{tabular}\\
\begin{tabular}{c} \input{pictures/det27.tex} \\$13\_13\_2_{||}2_{||}$\end{tabular}& \begin{tabular}{c} \input{pictures/ball27_1.tex} \\27A\end{tabular}& \begin{tabular}{c} \input{pictures/ball27_2.tex} \\27B\end{tabular}\\
\begin{tabular}{c} \input{pictures/det28.tex} \\$13\_13\_2_{||}4$\end{tabular}& \begin{tabular}{c} \input{pictures/ball28_1.tex} \\28A\end{tabular}& \begin{tabular}{c} \input{pictures/ball28_2.tex} \\28B\end{tabular}\\
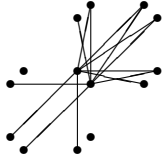
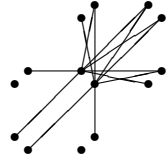
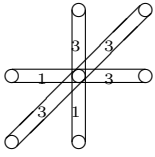
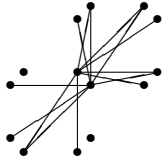
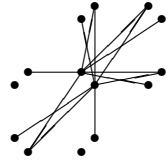
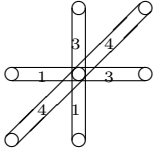
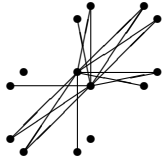
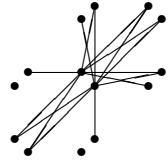
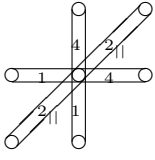
\begin{tabular}{c} \input{pictures/det29.tex} \\$13\_13\_33$\end{tabular}& \begin{tabular}{c} \input{pictures/ball29_1.tex} \\29A\end{tabular}& \begin{tabular}{c} \input{pictures/ball29_2.tex} \\29B\end{tabular}\\
\begin{tabular}{c} \input{pictures/det30.tex} \\$13\_13\_44$\end{tabular}& \begin{tabular}{c} \input{pictures/ball30_1.tex} \\30A\end{tabular}& \begin{tabular}{c} \input{pictures/ball30_2.tex} \\30B\end{tabular}\\
\begin{tabular}{c} \input{pictures/det31.tex} \\$14\_14\_2_{||}2_{||}$\end{tabular}& \begin{tabular}{c} \input{pictures/ball31_1.tex} \\31\end{tabular}\\
\begin{tabular}{c} \input{pictures/det32.tex} \\$14\_14\_2_{||}4$\end{tabular}& \begin{tabular}{c} \input{pictures/ball32_1.tex} \\32\end{tabular}\\
\begin{tabular}{c} \input{pictures/det33.tex} \\$14\_14\_33$\end{tabular}& \begin{tabular}{c} \input{pictures/ball33_1.tex} \\33\end{tabular}\\
\begin{tabular}{c} \input{pictures/det34.tex} \\$14\_14\_44$\end{tabular}& \begin{tabular}{c} \input{pictures/ball34_1.tex} \\34\end{tabular}\\
\begin{tabular}{c} \input{pictures/det35.tex} \\$2_{||}2_{||}\_2_{||}2_{||}\_33$\end{tabular}& \begin{tabular}{c} \input{pictures/ball35_1.tex} \\35\end{tabular}\\
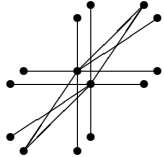
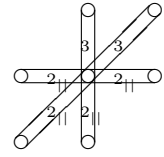
\begin{tabular}{c} \input{pictures/det36.tex} \\$2_{||}2_{||}\_2_{||}3\_2_{||}3$\end{tabular}& \begin{tabular}{c} \input{pictures/ball36_1.tex} \\36\end{tabular}\\
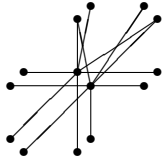
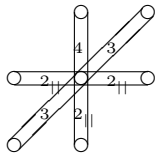
\begin{tabular}{c} \input{pictures/det37.tex} \\$2_{||}2_{||}\_2_{||}4\_33$\end{tabular}& \begin{tabular}{c} \input{pictures/ball37_1.tex} \\37\end{tabular}\\
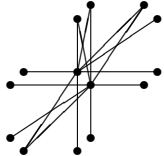
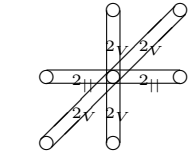
\begin{tabular}{c} \input{pictures/det38.tex} \\$2_{||}2_{||}\_2_V2_V\_2_V2_V$\end{tabular}& \begin{tabular}{c} \input{pictures/ball38_1.tex} \\38A\end{tabular}& \begin{tabular}{c} \input{pictures/ball38_2.tex} \\38B\end{tabular}\\
\begin{tabular}{c} \input{pictures/det39.tex} \\$2_{||}2_{||}\_33\_33$\end{tabular}& \begin{tabular}{c} \input{pictures/ball39_1.tex} \\39A\end{tabular}& \begin{tabular}{c} \input{pictures/ball39_2.tex} \\39B\end{tabular}\\
\begin{tabular}{c} \input{pictures/det40.tex} \\$2_{||}2_{||}\_33\_44$\end{tabular}& \begin{tabular}{c} \input{pictures/ball40_1.tex} \\40\end{tabular}\\
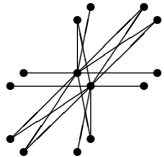
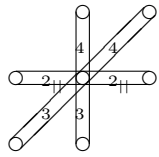
\begin{tabular}{c} \input{pictures/det41.tex} \\$2_{||}2_{||}\_34\_34$\end{tabular}& \begin{tabular}{c} \input{pictures/ball41_1.tex} \\41\end{tabular}\\
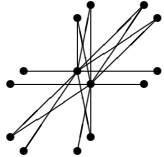
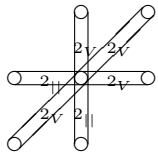
\begin{tabular}{c} \input{pictures/det42.tex} \\$2_{||}2_V\_2_{||}2_V\_2_V2_V$\end{tabular}& \begin{tabular}{c} \input{pictures/ball42_1.tex} \\42A\end{tabular}& \begin{tabular}{c} \input{pictures/ball42_2.tex} \\42B\end{tabular}\\
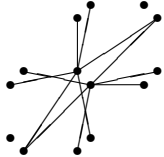
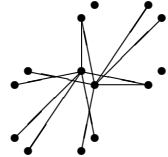
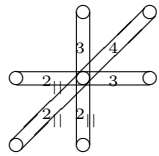
\begin{tabular}{c} \input{pictures/det43.tex} \\$2_{||}3\_2_{||}3\_2_{||}4$\end{tabular}& \begin{tabular}{c} \input{pictures/ball43_1.tex} \\43\end{tabular}\\
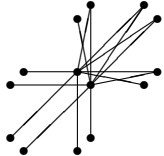
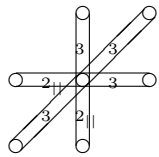
\begin{tabular}{c} \input{pictures/det44.tex} \\$2_{||}3\_2_{||}3\_33$\end{tabular}& \begin{tabular}{c} \input{pictures/ball44_1.tex} \\44\end{tabular}\\
\begin{tabular}{c} \input{pictures/det45.tex} \\$2_{||}3\_2_{||}3\_44$\end{tabular}& \begin{tabular}{c} \input{pictures/ball45_1.tex} \\45\end{tabular}\\
\begin{tabular}{c} \input{pictures/det46.tex} \\$2_{||}4\_2_{||}4\_33$\end{tabular}& \begin{tabular}{c} \input{pictures/ball46_1.tex} \\46\end{tabular}\\
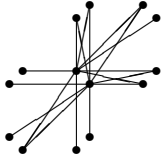
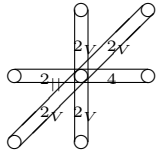
\begin{tabular}{c} \input{pictures/det47.tex} \\$2_{||}4\_2_V2_V\_2_V2_V$\end{tabular}& \begin{tabular}{c} \input{pictures/ball47_1.tex} \\47A\end{tabular}& \begin{tabular}{c} \input{pictures/ball47_2.tex} \\47B\end{tabular}\\
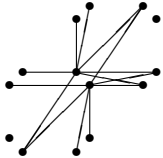
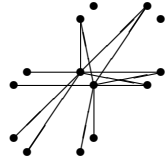
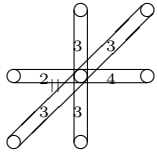
\begin{tabular}{c} \input{pictures/det48.tex} \\$2_{||}4\_33\_33$\end{tabular}& \begin{tabular}{c} \input{pictures/ball48_1.tex} \\48A\end{tabular}& \begin{tabular}{c} \input{pictures/ball48_2.tex} \\48B\end{tabular}\\
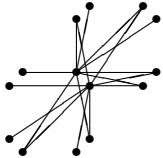
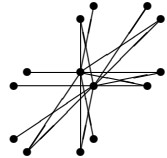
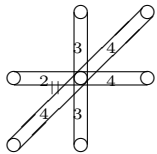
\begin{tabular}{c} \input{pictures/det49.tex} \\$2_{||}4\_33\_44$\end{tabular}& \begin{tabular}{c} \input{pictures/ball49_1.tex} \\49\end{tabular}\\
\begin{tabular}{c} \input{pictures/det50.tex} \\$2_{||}4\_34\_34$\end{tabular}& \begin{tabular}{c} \input{pictures/ball50_1.tex} \\50\end{tabular}\\
\begin{tabular}{c} \input{pictures/det51.tex} \\$2_V2_V\_2_V2_V\_33$\end{tabular}& \begin{tabular}{c} \input{pictures/ball51_1.tex} \\51A\end{tabular}& \begin{tabular}{c} \input{pictures/ball51_2.tex} \\51B\end{tabular}\\
\begin{tabular}{c} \input{pictures/det52.tex} \\$2_V2_V\_2_V2_V\_44$\end{tabular}& \begin{tabular}{c} \input{pictures/ball52_1.tex} \\52A\end{tabular}& \begin{tabular}{c} \input{pictures/ball52_2.tex} \\52B\end{tabular}\\
\begin{tabular}{c} \input{pictures/det53.tex} \\$2_V2_V\_2_V3\_2_V3$\end{tabular}& \begin{tabular}{c} \input{pictures/ball53_1.tex} \\53A\end{tabular}& \begin{tabular}{c} \input{pictures/ball53_2.tex} \\53B\end{tabular}& \begin{tabular}{c} \input{pictures/ball53_3.tex} \\53C\end{tabular}\\
\begin{tabular}{c} \input{pictures/det54.tex} \\$2_V2_V\_2_V4\_2_V4$\end{tabular}& \begin{tabular}{c} \input{pictures/ball54_1.tex} \\54A\end{tabular}& \begin{tabular}{c} \input{pictures/ball54_2.tex} \\54B\end{tabular}\\
\begin{tabular}{c} \input{pictures/det55.tex} \\$33\_33\_33$\end{tabular}& \begin{tabular}{c} \input{pictures/ball55_1.tex} \\55A\end{tabular}& \begin{tabular}{c} \input{pictures/ball55_2.tex} \\55B\end{tabular}\\
\begin{tabular}{c} \input{pictures/det56.tex} \\$33\_33\_44$\end{tabular}& \begin{tabular}{c} \input{pictures/ball56_1.tex} \\56A\end{tabular}& \begin{tabular}{c} \input{pictures/ball56_2.tex} \\56B\end{tabular}\\
\begin{tabular}{c} \input{pictures/det57.tex} \\$33\_34\_34$\end{tabular}& \begin{tabular}{c} \input{pictures/ball57_1.tex} \\57\end{tabular}\\
\begin{tabular}{c} \input{pictures/det58.tex} \\$33\_44\_44$\end{tabular}& \begin{tabular}{c} \input{pictures/ball58_1.tex} \\58\end{tabular}\\
\begin{tabular}{c} \input{pictures/det59.tex} \\$34\_34\_44$\end{tabular}& \begin{tabular}{c} \input{pictures/ball59_1.tex} \\59\end{tabular}\\
\end{longtable}

A list of 2872 saturated realizations generated using Algorithm 1 and thinned out using Algorithm 2
is given below in Table 4.
Saturated realizations are defined by triples $ H, L, X $ (see above):
group $H\in \textbf{H}$ is given in the fourth column,
subgroup $L$ of index $ 2 $ in $ H _ {(0,0,0)} $ - in the fifth column,
a subset of $X$ of elements of the group $H$ is in the seventh column.
In addition, in the sixth column is given an element $m$, such that
$L\cup mL = H_{(0,0,0)}$. The realizaion \textnumero \, is given in the third column.
Saturated realizations are sorted by vertex neighborhood extensions given in Table 3.
(\textnumero of vertex neighborhood extension is given in the first column of Table 4).
The set of realizations with given vertex neighborhood extension, in turn, 
are divided into classes,
defined by orders of balls of radii 1, 2, ..., 10 of graph of a realization
(we call such a set of orders of balls of radii 1, 2, ..., 10 by {\it growth} and give it in the second column of Table 4).
In this subdivision, classes are lexicographically sorted by growth increasing.
In Table 4, along with saturated realizations, we give non-saturated realizations by
\textnumero of saturated realizations with an asterisk in the third column (details see before Corollary 1 below).

\hspace*{97mm}{\small \mbox{T\ a\ b\ l\ e\ \  4}}\\
\centerline{\small\bf 2872 saturated and 2701 non-saturated maximal realizations}\\
 \centerline{\small\bf of symmetrical 2-extensions of the grid $\Lambda^{3}$ of class I}
{\tiny
\begin{longtable}{p{0.01\linewidth}p{0.01\linewidth}llp{0.12\linewidth}p{0.03\linewidth}p{0.55\linewidth}}
Nbr. & gr & \textnumero & $H_i$ & $L$ & $m$ & $X$ \\
\hline\endhead
\\
\multicolumn{7}{l}{\normalsize 1A}\\
& \multicolumn{6}{l}{[ 3, 6, 12, 20, 32, 50, 65, 82, 105, 132 ]}\\
& & \multicolumn{5}{l}{20*, }\\
& \multicolumn{6}{l}{[ 4, 9, 19, 35, 52, 72, 100, 131, 163, 201 ]}\\
& & 1 & $H_{778}$	& $\langle m_z r_x, r_y r_x \rangle$	& $r_z^2 r_x$	&
 $ m_x t_x^{-1}, r_z^{-1} t_y^{-1}, m_x r_y^{-1} r_z^{-1} t_z^{-1} $\\ \hline& \multicolumn{6}{l}{[ 4, 10, 21, 36, 54, 78, 106, 136, 173, 214 ]}\\
& & 2 & $H_{657}$	& $\langle m_z r_x \rangle$	& $r_z^2 r_x$	&
 $ m_x t_x^{-1}, r_x^2 t_y^{-1}, m_z r_x^{-1} t_z^{-1} $\\
\hline
 \end{longtable}}

\begin{theorem}
Saturated realizations of symmetrical 2-extensions of the grid $\Lambda^{3}$ of class I up to equivalence
are exhausted by $2872$ pairwise nonequivalent saturated realizations given in Table 4.
\end{theorem}

Note that listing of all, up to equivalence, realizations of the symmetrical 2-extensions of the grid $\Lambda^{3}$
is reduced to listing of all, up to equivalence, saturated realizations of symmetrical 2-extensions of $\Lambda^{3}$.
Indeed, it is obvious that every non-saturated realization of a symmetrical 2-extension of $\Lambda^{3}$
is obtained from a uniquely defined saturated realization by removing of an edge in each block.
All maximal non-saturated realizations obtained in this way are given in Table 4 by \textnumero \, 
of saturated maximal realizaions taken with an asterisk.

\begin{corollary} 
Non-saturated realizations of symmetrical 2-extensions of the grid $\Lambda^{3}$ of class I up to equivalence
are exhausted by $2701$ pairwise nonequivalent non-saturated realizations given in Table 4.
\end{corollary}
\begin{proof}
Using a computer, it is easily verified that when removing edges inside blocks of $2872$ realizaions given in Table 4,
171 their graphs become disconnected,
and the remaining 2701 realizaions (see \textnumero \,  with an asterisk in Table 4) give all, up to equivalence,
non-saturated realizations of symmetrical 2-extensions of the grid $\Lambda^{3}$ of class I.
\end{proof}

According to Theorem 1 and Corollary 1, all realizations given in Table 4 are pairwise nonequivalent.
However, among the graphs of these realizaions, isomorphic ones are found.
With GAP, we built
partition of the set of graphs of realizaions from Table 4 into classes of isomorphic graphs
(for details, see the proof of Corollary \ref{cor2} below).
In the following Table 5, we give all the non-singleton classes of this partition.

\hspace*{97mm}{\small \mbox{T\ a\ b\ l\ e\ \  5}}\\
\centerline{\small\bf Non-singleton classes of ishomorphic graphs of realizaions}\\
 \centerline{\small\bf of symmetrical 2-extensions of $\Lambda^{3}$ of class I}

{\tiny
\begin{multicols}{4}
{
\noindent
1) 1, 45*\\
2) 7, 8, 76*\\
3) 9, 10\\
4) 16, 17\\
5) 18, 19\\
6) 26, 64*\\
7) 30, 97*\\
8) 36, 198*, 199*, 337*\\
9) 37, 38, 202*, 203*, 204*, 205*, 206*, 207*, 1119*, 1120*, 1484*\\
10) 39, 1121*\\
11) 41, 279\\
12) 42, 280\\
13) 43, 44, 214*, 215*, 343*\\
14) 46, 123*\\
15) 47, 122*\\
16) 54, 293\\
17) 55, 294\\
18) 57, 348*, 352*\\
19) 59, 144*\\
20) 64, 65\\
21) 110, 250*\\
22) 122, 123\\
23) 198, 199\\
24) 202, 203\\
25) 206, 207\\
26) 216, 343\\
27) 229, 348\\
28) 233, 352\\
29) 286, 216*\\
30) 296, 1123*, 1124*\\
31) 299, 229*\\
32) 301, 233*\\
33) 324, 363*\\
34) 332, 366*\\
35) 427, 2114*\\
36) 428, 503*, 658*\\
37) 431, 667*, 1489*\\
38) 439, 1422*\\
39) 440, 1421*\\
40) 455, 1440*\\
41) 456, 1439*\\
42) 651, 652\\
43) 653, 654\\
44) 664, 1489\\
45) 672, 1421\\
46) 673, 1422\\
47) 712, 1439\\
48) 713, 1440\\
49) 1091, 1094\\
50) 1092, 1093\\
51) 1282, 1412*\\
52) 1283, 1414*\\
53) 1288, 664*, 1423*\\
54) 1289, 1420*\\
55) 1290, 1424*\\
56) 1294, 672*\\
57) 1295, 673*\\
58) 1307, 1437*\\
59) 1313, 712*\\
60) 1314, 713*\\
61) 1716, 1717\\
62) 1720, 1721\\
63) 1722, 1724\\
64) 1725, 1726\\
65) 1727, 1728\\
66) 1756, 1757\\
67) 1759, 1760\\
68) 1762, 1763\\
69) 1857, 1858\\
70) 1919, 1920\\
71) 1928, 1929\\
72) 1941, 1943\\
73) 1944, 1945\\
74) 1962, 1963\\
75) 1990, 1991\\
76) 1998, 1999\\
77) 2002, 2006\\
78) 2004, 2005\\
79) 2043, 2047\\
80) 2050, 2051\\
81) 2053, 2054\\
82) 2055, 2056\\
83) 2058, 2059\\
84) 2071, 2085\\
85) 2072, 2086\\
86) 2075, 2084\\
87) 2080, 2082\\
88) 2089, 2090\\
89) 2091, 2092\\
90) 2093, 2094\\
91) 2109, 2110\\
92) 2346, 2349\\
93) 2347, 2348\\
94) 2353, 2355\\
95) 2357, 2358\\
96) 2362, 2365\\
97) 2363, 2364\\
98) 2368, 2369\\
99) 2371, 2372\\
100) 2374, 2375\\
101) 2376, 2378\\
102) 2492, 2495\\
103) 2552, 2553\\
104) 2651, 2659\\
105) 2652, 2660\\
106) 2653, 2658\\
107) 2654, 2655\\
108) 2656, 2657\\
109) 2702, 2703\\
110) 2707, 2708\\
111) 2713, 2714\\
112) 2718, 2720\\
113) 2724, 2726\\
114) 63*, 64*, 65*\\
115) 89*, 314*\\
116) 102*, 105*\\
117) 116*, 1121*\\
118) 117*, 339*\\
119) 121*, 219*\\
120) 130*, 1127*\\
121) 141*, 1125*\\
122) 146*, 353*\\
123) 167*, 245*\\
124) 175*, 1139*\\
125) 209*, 210*\\
126) 211*, 212*\\
127) 228*, 350*\\
128) 234*, 235*\\
129) 328*, 329*\\
130) 360*, 361*\\
131) 651*, 652*, 653*, 654*, 2114*\\
132) 660*, 2115*\\
133) 662*, 2116*\\
134) 817*, 818*\\
135) 819*, 820*\\
136) 821*, 2117*\\
137) 825*, 2118*\\
138) 1091*, 1094*\\
139) 1092*, 1093*\\
140) 1133*, 1134*\\
141) 1136*, 1137*\\
142) 1143*, 1144*\\
143) 1147*, 1148*\\
144) 1149*, 1150*\\
145) 1152*, 2521*\\
146) 1163*, 1622*\\
147) 1716*, 1717*\\
148) 1720*, 1721*\\
149) 1722*, 1724*\\
150) 1725*, 1726*\\
151) 1727*, 1728*\\
152) 1756*, 1757*\\
153) 1759*, 1760*\\
154) 1762*, 1763*\\
155) 1857*, 1858*\\
156) 1919*, 1920*\\
157) 1928*, 1929*\\
158) 1941*, 1943*\\
159) 1944*, 1945*\\
160) 1962*, 1963*\\
161) 1990*, 1991*\\
162) 1998*, 1999*\\
163) 2002*, 2006*\\
164) 2004*, 2005*\\
165) 2043*, 2047*\\
166) 2050*, 2051*\\
167) 2053*, 2054*\\
168) 2055*, 2056*\\
169) 2058*, 2059*\\
170) 2071*, 2085*\\
171) 2072*, 2086*\\
172) 2075*, 2084*\\
173) 2080*, 2082*\\
174) 2089*, 2090*\\
175) 2091*, 2092*\\
176) 2093*, 2094*\\
177) 2109*, 2110*\\
178) 2346*, 2349*\\
179) 2347*, 2348*\\
180) 2353*, 2355*\\
181) 2357*, 2358*\\
182) 2362*, 2365*\\
183) 2363*, 2364*\\
184) 2368*, 2369*\\
185) 2371*, 2372*\\
186) 2374*, 2375*\\
187) 2376*, 2378*\\
188) 2492*, 2495*\\
189) 2552*, 2553*\\
190) 2651*, 2659*\\
191) 2652*, 2660*\\
192) 2653*, 2658*\\
193) 2654*, 2655*\\
194) 2656*, 2657*\\
195) 2702*, 2703*\\
196) 2707*, 2708*\\
197) 2713*, 2714*\\
198) 2718*, 2720*\\
199) 2724*, 2726*\\ 
} \end{multicols} }

\begin{corollary}\label{cor2} 
$(1)$ Up to isomorphism, there are $2792$ graphs of saturated realizations of symmetrical 2-extensions of the grid $\Lambda^{3}$ of class I.\\ 
\\ 
$(2)$ Up to isomorphism, there are $2594$ graphs of non-saturated realizations 
of symmetrical 2-extensions of the grid $\Lambda^{3}$ of class I.\\  
$(3)$ Up to isomorphism, there are $5350$ graphs of realizations 
of symmetrical 2-extensions of the grid $\Lambda^{3}$ of class I.\\ 
\end{corollary}

\begin{proof}
Using GAP for each of the graphs of 5573 (2872 saturated and 2701 non-saturated) realizaions 
of symmetrical 2-extensions of $\Lambda^{3}$ of class I
a subgraph $B$ was generated by
a set of vertices, which are at a distance of $\leq 4$ from some arbitrary vertex (i.e. $B$ is a ball of radius 4).
In the obtained set of 5573 finite graphs, balls are isomorphic if and only if
they correspond to realizaions which are in the same line of Table 5.
After that, each isomorphism $\varphi_b$ between the balls $B_1$ and $B_2$ we continued to the isomorphism $\varphi$
of whole graphs of the corresponding realizaions $R_1 = (\Gamma_1, G_1, $ $ \sigma_1, \varphi_1)$ and
$R_2 = (\Gamma_2, G_2, \sigma_2, \varphi_2)$ as follows.

Let a realization $R_1$ satisfy the condition of $[p_{x1}, p_{y1}, p_{z1}]$ - periodicity,
and $R_2$ -- the condition of $[p_{x2}, p_{y2}, p_{z2}]$ - periodicity
(according to \cite{IV2} a realizaion $R = (\Gamma, G, \sigma, \varphi)$ of a symmetrical 2-extension of $\Lambda^{3}$
{\it satisfies the condition of $[p_x, p_y, p_z]$-periodicity}, where $p_x, p_y, p_z$ are positive integers, if
there exist $g_1, g_2, g_3 \in G$ such that $[g_i, g_j] = 1$ for $i\neq j$ and
$\varphi g_1^{\sigma} \varphi^{-1} = t_x^{p_x}$, $\varphi g_2^{\sigma} \varphi^{-1} = t_y^{p_y}$, 
$\varphi g_3^{\sigma} \varphi^{-1} = t_z^{p_y}$).
We identify the set of vertices of $\Gamma_1$ with the set
$\{(x,y,z,w): x,y,z\in \mathbb Z, w\in\{0,1\}\}$, so that $\sigma_1=\{\{(x,y,z,0),(x,y,z,1)\}: x,y,z\in \mathbb Z\}$
and $\{(x_1,y_1,z_1,w_1),(x_2,y_2,z_2,w_2)\} \in \mathrm E(\Gamma_1) 
\Leftrightarrow \{(x_1+p_{x1},y_1,z_1,w_1),(x_2+p_{x1},y_2,z_2,w_2)\} \in \mathrm E(\Gamma_1)
\Leftrightarrow \{(x_1,y_1+p_{y1},z_1,w_1)$, $(x_2,y_2+p_{y1},z_2,w_2)\} \in \mathrm E(\Gamma_1)
\Leftrightarrow \{(x_1,y_1,z_1+p_{z1},w_1),(x_2,y_2,z_2+p_{z1},w_2)\} \in \mathrm E(\Gamma_1)$.
Similarly, we identify the set of vertices of $\Gamma_2$ with the set
$\{(x,y,z,w): x,y,z\in \mathbb Z, w\in\{0,1\}\}$.

We select positive integers $p_x, p_y, p_z$, so that
the realization $R_1$ satisfies the conditions of $[p_x, p_y, p_z]$-periodicity,
$p_x | p_{x1}, p_y | p_{y1}, p_z | p_{z1}$ and
$(p_x, 0, 0) M$, $(0, p_y, 0) M$, $(0, 0, p_z) M \in \langle(p_{x2}, 0,0), (0,p_{y2},0), (0,0,p_{z2})\rangle$
(angle brackets mean generation in the additive group of row-vectors), where
the $3\times 3$-matrix $M$ is obtained from the $3\times 4$ matrix
$$\left(\begin{array}{c}
(p_x,0,0,0)\varphi_b/p_x\\
(0,p_y,0,0)\varphi_b/p_y\\
(0,0,p_z,0)\varphi_b/p_z
\end{array}\right)$$
by deleting of the last column.
We need to ensure that the extension of the fragment $[0...p_x-1]\times[0...p_y-1]\times[0...p_z-1]$
 is inside of the ball $B_1$. To choose $p_x, p_y, p_z$ in such way for some
 pairs of realizations $R_1$ and $R_2$ we had to take the radius of balls $B_1$ and $B_2$ greater than 4.

The image of an arbitrary vertex $u$ of the extension $\Gamma_1$ is now defined by
 $u\varphi=$\linebreak $=u t^{-1} \varphi_b t^M$, where the shift $t^{-1}=t_x^{p_x n_1} t_y^{p_y n_2} t_z^{p_z n_3}$ 
 maps $u$ into the extension of fragment $[0...p_x-1]\times[0...p_y-1]\times[0...p_z-1]$
 ($n_1, n_2, n_3$ are suitable positive integers).
\end{proof}


\Addresses
\end{document}